\documentclass{article}
\usepackage{amssymb}

\usepackage{graphicx}
\usepackage{amsmath}

\newtheorem{theorem}{Theorem}

\newtheorem{remark}[theorem]{Remark}

\input{tcilatex}

\begin{document}

\title{From Special Lagrangian to Hermitian-Yang-Mills\\
via Fourier-Mukai Transform}
\author{Naichung Conan Leung, Shing-Tung Yau, and Eric Zaslow}
\maketitle

\begin{abstract}
We exhibit a transformation taking special Lagrangian submanifolds of a
Calabi-Yau together with local systems to vector bundles over the mirror
manifold with connections obeying deformed Hermitian-Yang-Mills equations.
That is, the transformation relates supersymmetric A- and B-cycles. In this
paper, we assume that the mirror pair are dual torus fibrations with flat
tori and that the A-cycle is a section.

We also show that this transformation preserves the (holomorphic)
Chern-Simons functional for all connections. Furthermore, on corresponding
moduli spaces of supersymmetric cycles it identifies the graded tangent
spaces and the holomorphic m-forms. In particular, we verify Vafa's mirror
conjecture with bundles in this special case.
\end{abstract}

\tableofcontents

\section{Introduction}

In this note we argue that a real version of the Fourier-Mukai transform
would carry supersymmetric A-cycles to B-cycles. Roughly, special Lagrangian
submanifolds (plus local systems) will be mapped to holomorphic submanifolds
(plus bundles over them). The approach is similar to the one in \cite
{polarin}, though our emphasis is geometric, as we focus on the differential
equations defining D-branes.

Our goal is to provide the basis for a geometric functor relating the
categories of D-branes on opposite mirror sides: the derived category of
coherent sheaves on one hand, and Fukaya's category of Lagrangian
submanifolds with local systems on the other. (In fact, we consider special
Lagrangians as objects.\footnote{%
One can't yet say which formulation of these categories will best fit the
physics. The Donaldson-Uhlenbeck-Yau theorem relates stable bundles to
solutions of the Hermitian-Yang-Mills equations. Analogously, one hopes that
Lagrangians are equivalent to special Lagrangians up to Hamiltonian
deformations. We are a long way away from proving such theorems, however.
The recent preprint of R. Thomas investigates these
two issues as moment map problems related by mirror symmetry.
})
At this point, providing a physical interpretation of the derived category
is premature -- even though recent discussions of the role of the
brane-anti-brane tachyon offer glimpses -- so we content ourselves with
thinking of D-branes as vector bundles over holomorphic submanifolds, where
``submanifold'' may mean the entire space.

We assume, following \cite{syz}, that the mirror manifold has a description
as a dual torus fibration, blithely ignoring singular fibers for now. The
procedure is now quite simple \cite{polarin}\cite{Donagi}\cite{fmw}: a point
on a torus determines a line bundle on a dual torus. A section of a torus
fibration then determines a family of line bundles on the mirror side. These
can fit together to define a bundle. We will show how, in this idealized
situation, the differential equations on the two sides are related under
this transformation. We also show that the Chern-Simons action of an A-cycle
equals the holomorphic Chern-Simons action of its transform, even off-shell.
Vafa's version of mirror symmetry with bundles is then verified in this
setting.

It will be interesting, though perhaps quite formidable, to generalize this
procedure by relaxing some of our assumptions and extending the setting to
include singular fibers and more general objects in the derived category.

\textbf{Acknowledgements:} We gratefully acknowledge helpful conversations
with Richard Thomas, who
has also obtained similar structures for moduli spaces
of supersymmetric cycles
in his recent preprint (which he has kindly shared with us).
N.C. Leung is supported by NSF grant DMS-9803616
and S.-T. Yau is supported by NSF grant DMS-9803347.

\section{Supersymmetric A- and B-cycles}

The authors of \cite{mmms} consider the supersymmetric p-brane action and
determine the conditions for preserving supersymmetry (BPS).\footnote{%
In this section, we follow what has become standard notation. Let us remark,
however, that the IIA string theory, on which one normally considers the
A-model, naturally has even-dimensional (hence type-B) branes, so-named
because their compositions depend on the complex structure. IIB string
theory has odd-dimensional (type-A) branes, whose composition depends only
on the symplectic structure.} They show that there are two kinds of
supersymmetric cycles $(C,L)$ on a Calabi-Yau threefold $M$ where $C$ is a
(possibly singular and with multiplicity) submanifold of $M$ and $L$ is a
complex line bundle over $C$ together with a $U\left( 1\right) $ connection $%
D_{A}$. Let us denote the K\"ahler form (resp. holomorphic volume form) on
the Calabi-Yau threefold by $\omega $ (resp. $\Omega $).

The type-A supersymmetric cycle is when $C$ is a special Lagrangian
submanifold of $M$ and the curvature $F_{A}$ of $D_{A}$ vanishes, 
\begin{equation*}
F_{A}=0,
\end{equation*}
namely $D_{A}$ is a unitary flat connection. In the presence of a background 
$B$-field (an element of $H^{2}(M,\mathbf{{R}/{Z})}$), $F_{A}$ should be
replaced by $F_{A}-B,$ where $B$ is understood to be pulled back to the
submanifold. We take $B=0$ in this paper. Recall that a Lagrangian
submanifold $C$ is called ``special'' if when restricting to $C$ we have 
\begin{equation*}
\func{Im}\Omega =\tan \theta \func{Re}\Omega ,
\end{equation*}
for some constant $\theta .$ Or equivalently, $\func{Im}e^{i\theta }\Omega
=0 $.

The type-B cycle is when $C$ is a complex submanifold of $M$ of dimension $n$
and the curvature two form $F_{A}$ of $D_{A}$ satisfies following
conditions: 
\begin{eqnarray*}
F_{A}^{0,2} &=&0, \\
\func{Im}e^{i\theta }\left( \omega +F_{A}\right) ^{n} &=&0.
\end{eqnarray*}
The first equation says that the (0,1) component of the connection
determines a holomorphic structure on $L$. The second equation is called the
deformed Hermitian-Yang-Mills equation and it is equivalent to the following
equation, 
\begin{equation*}
\func{Im}\left( \omega +F_{A}\right) ^{n}=\tan \theta \func{Re}\left( \omega
+F_{A}\right) ^{n}.
\end{equation*}
For example when $C$ is the whole Calabi-Yau manifold $M$ of dimension three
then the second equation says $F\wedge \omega ^{2}/2-F^{3}/6=\tan \theta %
\left[ \omega ^{3}/6-\left( F^{2}/2\right) \wedge \omega \right] $.

\section{Fourier-Mukai Transform of A- and B-cycles}

In this section we explain the Fourier-Mukai transform of supersymmetric
cycles. The gist of the story is that, assuming mirror pairs are mirror
torus fibrations, each point of a Lagrangian submanifold lies in some fiber
-- hence defines a bundle over the dual fiber. When done in families and
with connections, we get a bundle with connection on the mirror, and the
differential equations defining A-cycles map to those which define B-cycles
on the mirror. Recall that the base of the fibration itself -- the zero
graph -- should be dual to the six-brane with zero connection.
Multi-sections are dual to higher-rank bundles, and are discussed in section
3.2. Other cases appear in 3.3.

We assume that the $m$ dimensional Calabi-Yau mirror pair $M$ and $W$ have
dual torus fibrations. To avoid the difficulties of singular fibers and
unknown Calabi-Yau metrics, we will only consider a neighborhood of a smooth
special Lagrangian torus and also assume the K\"{a}hler potential $\phi $ on 
$M$ to be $T^{m}$-invariant (see for example p.20 of \cite{hitchin}). This
is the semi-flat assumption of \cite{syz}. Notice that the Lagrangian
fibrations on $M$ and $W$ are in fact special.

Therefore, let $\phi \left( x^{j},y^{j}\right) =\phi \left( x^{j}\right) $. (%
$y$ is the coordinate for the fiber and $x$ for the base $B$ of the
fibration on $M$. The holomorphic coordinates on $M$ are $z^{j}=x^{j}+iy^{j}$%
's.) As studied by Calabi, the Ricci tensor vanishes and $\Omega
=dz^{1}\wedge ....\wedge dz^{m}$ is covariant constant if and only if $\phi $
satisfies a real Monge-Amp\`{e}re equation 
\begin{equation*}
\det \frac{\partial ^{2}\phi }{\partial x^{i}\partial x^{j}}=const.
\end{equation*}
The Ricci-flat K\"{a}hler metric and form are 
\begin{eqnarray*}
g &=&\sum_{i,j}\frac{\partial ^{2}\phi }{\partial x^{i}\partial x^{j}}\left(
dx^{i}dx^{j}+dy^{i}dy^{j}\right) , \\
\omega &=&\frac{i}{2}\sum_{i,j}\frac{\partial ^{2}\phi }{\partial
x^{i}\partial x^{j}}dz^{i}\wedge d\overline{z}^{j}\text{,}
\end{eqnarray*}
(henceforth we sum over repeated indices). Notice that $\Omega \wedge \bar{%
\Omega}$ is a constant mulitple of $\omega ^{m}$ and it is direct
consequence of the real Monge-Amp\`{e}re equation.

Also note from the form of the metric $g$ that $M$ is locally isometric to
the tangent bundle of $B$ with its metric induced from the metric $\sum_{i,j}%
\frac{\partial ^{2}\phi }{\partial x^{i}\partial x^{j}}dx^{i}dx^{j}$ on $B$.
If we use the metric on $B$ to identify its tangent bundle with its
cotangent bundle, then the above symplectic form $\omega $ is just the
canonical symplectic form $dp\wedge dq$ on the cotangent bundle.

We can view the universal cover of $M$ either as $TB$ with the standard
complex structure, or as $T^{\ast }B$ with the standard symplectic
structure. A solution of the real Monge-Amp\`{e}re equation is used to
determine the symplectic structure in the former case and to determine the
metric structure, and therefore the complex structure in the latter case.

\subsection{Transformation of a section}

We will construct the transform for a special Lagrangian exhibited as a
section of the fibration, i.e.\negthinspace\ a graph over the base.

Recall that a section of $T^{\ast }B$ is Lagrangian with respect to the
standard symplectic form if and only if it is a closed one form, and hence
locally exact. Therefore (or by calculation), a graph $y(x)$ in $M$ is
Lagrangian with respect to $\omega $ if and only if $\frac{\partial }{%
\partial x^{j}}(y^{l}\phi _{lk})=\frac{\partial }{\partial x^{k}}(\phi
_{lj}y^{l}),$ where $\phi _{ij}=\frac{\partial ^{2}\phi }{\partial
x^{i}\partial x^{j}}$, from which we get 
\begin{equation*}
y^{j}=\phi ^{jk}\frac{\partial f}{\partial x^{k}}
\end{equation*}
for some function $f$ (locally), where $\phi ^{jk}$ is the inverse matrix of 
$\phi _{jk}$.

Now $dz^{j}=dx^{j}+idy^{j}$ and on $C$ we have $dy^{j}=\phi ^{jl}\left( 
\frac{\partial ^{2}f}{\partial x^{l}\partial x^{k}}-\phi ^{pq}\phi _{lkp}%
\frac{\partial f}{\partial x^{q}}\right) dx^{k}$. Therefore $dz^{j}=\left(
\delta _{jk}+i\phi ^{jl}\left( \frac{\partial ^{2}f}{\partial x^{l}\partial
x^{k}}-\phi ^{pq}\phi _{lkp}\frac{\partial f}{\partial x^{q}}\right) \right)
dx^{k}$ over $C$. Notice that if we write $g=\phi _{jk}dx^{j}dx^{k}$ as the
Riemannian metric on the base, then the Christoffel symbol for the
Levi-Civita connection is $\Gamma _{lk}^{q}=\phi ^{pq}\phi _{lkp}.$
Therefore $Hess\left( f\right) =\left( \frac{\partial ^{2}f}{\partial
x^{l}\partial x^{k}}-\phi ^{pq}\phi _{lkp}\frac{\partial f}{\partial x^{q}}%
\right) dx^{l}dx^{k}.$ Hence 
\begin{eqnarray*}
dz^{1}\wedge ...\wedge dz^{m}|_{C} &=&\det \left( I+ig^{-1}Hess\left(
f\right) \right) dx^{1}\wedge ...\wedge dx^{m} \\
&=&\det \left( g\right) ^{-1}\det \left( g+iHess\left( f\right) \right)
dx^{1}\wedge ...\wedge dx^{m},
\end{eqnarray*}
so the special Lagrangian condition (with phase) $\func{Im}\left(
dz^{1}\wedge ....\wedge dz^{m}\right) |_{C}=\tan \theta \cdot Re\left(
dz^{1}\wedge ....\wedge dz^{m}\right) |_{C}$ becomes 
\begin{equation*}
\func{Im}\det \left( g+iHess\left( f\right) \right) =\left( \tan \theta
\right) \func{Re}\det \left( g+iHess\left( f\right) \right) .
\end{equation*}

{}From these data, we want to construct a $U\left( 1\right) $ connection
over the mirror manifold $W$ which satisfies the deformed
Hermitian-Yang-Mills equation. The dual manifold $W$ is constructed by
replacing each torus fiber $T$ in $M$ by the dual torus $\widetilde{T}%
=Hom\left( T,S^{1}\right) $. If we write the dual coordinates to $%
y^{1},...,y^{m}$ as $\widetilde{y}_{1},...,\widetilde{y}_{m}$, then the dual
Riemannian metric on $W$ is obtained by taking the dual metric on each dual
torus fiber $\widetilde{T}$:

\begin{equation*}
\widetilde{g}=\sum_{i,j}\left( \phi _{ij}dx^{i}dx^{j}+\phi ^{ij}d\widetilde{y%
}_{i}d\widetilde{y}_{j}\right) .
\end{equation*}
We need to understand the complex structure and the symplectic structure on $%
W$ (see for example \cite{syz} and \cite{hitchin}). First we rewrite $%
\widetilde{g}$ as follows, 
\begin{equation*}
\widetilde{g}=\sum_{i,j}\phi ^{ij}\left( \left( \Sigma _{k}\phi
_{ik}dx^{k}\right) \left( \Sigma _{l}\phi _{jl}dx^{l}\right) +d\widetilde{y}%
_{i}d\widetilde{y}_{j}\right) .
\end{equation*}
Notice that $d\left( \Sigma _{k}\phi _{jk}dx^{k}\right) =0$ because $\phi
_{jkl}$ is symmetric with respect to interchanging the indexes. Therefore
there exist functions $\widetilde{x}_{j}=\widetilde{x}_{j}\left( x\right) $%
's such that $d\widetilde{x}_{j}=\Sigma _{k}\phi _{jk}dx^{k}$ locally --
then $\frac{\partial \widetilde{x}_{j}}{\partial x^{k}}=\phi_{jk}$ -- and we
obtain 
\begin{equation*}
\widetilde{g}=\sum_{i,j}\phi ^{ij}\left( d\widetilde{x}_{i}d\widetilde{x}%
_{j}+d\widetilde{y}_{i}d\widetilde{y}_{j}\right) .
\end{equation*}
So we can use $\widetilde{z}_{j}=\widetilde{x}_{j}+i\widetilde{y}_{j}$'s as
complex coordinates on $W$. It is easy to check that the corresponding
symplectic form is given by 
\begin{equation*}
\widetilde{\omega }= \frac{i}{2}\sum_{i,j}\phi ^{ij}d\widetilde{z}_{i}\wedge
d\overline{\widetilde{z}}_{j}.
\end{equation*}
Moreover the covariant constant holomorphic m-form on $W$ is given by 
\begin{equation*}
\widetilde{\Omega }=d\widetilde{z}_{1}\wedge ...\wedge d\widetilde{z}_{m}.
\end{equation*}
Again, as a direct consequence of $\phi $ being a solution of the real
Monge-Amp\`{e}re equation, $\widetilde{\Omega }\wedge \overline{\widetilde{%
\Omega }}$ is a constant multiple of $\widetilde{\omega }^{m}$.

\begin{remark}
The mirror manifold $W$ can be interpreted as the moduli space of special
Lagrangian tori together with flat $U(1)$ connections over them (see \cite
{syz}). It is because the dual torus parametrizes isomorphism classes of
flat $U(1)$ connections on the original torus. It can be checked directly
that the $L^{2}$ metric, i.e. the Weil-Petersson metric, on this moduli
space $W$ coincides with our $\widetilde{g}$ above.

In general, the relevent metric on the moduli space $W$ is given by a
two-point function computed via a path integral, wihch includes instanton
contributions from holomorphic disks bounding the special Lagrangian torus
fibers. However, for our local Calabi-Yau $M$ such holomorphic disks do not
exist. This is because $M$ is homotopic to any one of its fibers; but any
such holomorphic disk would define a non-trivial relative homology class.
Therefore our metric $\widetilde{g}$ coincides with the physical metric on
the moduli space $W$.
\end{remark}

\begin{remark}
We note the symmetry between $g$ (resp. $\omega $) and $\widetilde{g}$
(resp. $\widetilde{\omega }$). For one can write $\phi ^{ij}$ as the second
derivative of some function $\widetilde{\phi }$ with respect to the $%
\widetilde{x}_{j}$'s. Simply write $x^{j}=x^{j}\left( \widetilde{x}\right) $%
, then $\frac{\partial x^{j}}{\partial \widetilde{x}_{k}}=\phi ^{jk}=\frac{%
\partial x^{k}}{\partial \widetilde{x}_{j}}$ and therefore $x^{j}=\frac{%
\partial \Phi }{\partial \widetilde{x}_{j}}$ for some function, $\Phi ,$ and
it is easy to check that $\widetilde{\phi }=\Phi $.
\end{remark}

On each torus fiber, we have canonical isomorphisms $T =\mathrm{Hom}(%
\widetilde{T},S^{1})=\mathrm{Hom}( \pi _{1}( \widetilde{T}) ,S^{1}),$
therefore a point $y=\left( y^{1},...,y^{m}\right) $ in $T$ defines a flat
connection $D_{y}$ on its dual $\widetilde{T\text{.}}$ This is the real
Fourier-Mukai transform. Explicitly, we have 
\begin{equation*}
\begin{array}{cc}
g_{y}: & \widetilde{T}\rightarrow i\left( \mathbb{R}/\mathbb{Z}\right) =
S^{1} \\ 
& \widetilde{y}\mapsto i\sum_{j=1}^{m}y^{j}\widetilde{y}_{j},
\end{array}
\end{equation*}
and $D_{y}=d + A = d+idg_{y}=d+i\Sigma y^{j}d\widetilde{y}_{j}.$

In fact we get a torus family of one-forms, since $y$ (hence $A$) has $x$-
(or $\widetilde{x}$-) dependence. Namely, we obtain a $U\left( 1\right) $
connection on $W$, 
\begin{equation*}
D_{A}=d+i\sum_{j}y^{j}d\widetilde{y}_{j}.
\end{equation*}
Its curvature two form is given by, 
\begin{equation*}
F_{A}=dA=\sum_{k,j}i\frac{\partial y^{j}}{\partial \widetilde{x}_{k}}d%
\widetilde{x}_{k}\wedge d\widetilde{y}_{j}.
\end{equation*}
In particular 
\begin{equation*}
F_{A}^{2,0}=\frac{1}{2}\sum_{j,k}\left( \frac{\partial y^{k}}{\partial 
\widetilde{x}_{j}}-\frac{\partial y^{j}}{\partial \widetilde{x}_{k}}\right) d%
\widetilde{z}_{j}\wedge d{\widetilde{z}}_{k}.
\end{equation*}
Therefore, that $D_{A}$ is integrable, i.e. $F_{A}^{0,2}=0$, is equivalent
to the existence of $f=f\left( \widetilde{x}\right) $ such that $y^{j}=\frac{%
\partial f}{\partial \widetilde{x}_{j}}=\phi ^{jk}\frac{\partial f}{\partial
x^{j}}$ because of $d\widetilde{x}_{j}=\Sigma _{k}\phi _{jk}dx^{k}$. Namely,
the cycle $C\subset M$ must be Lagrangian. Now 
\begin{equation*}
\frac{\partial y^{j}}{\partial \widetilde{x}_{k}}=\frac{\partial ^{2}f}{%
\partial \widetilde{x}_{j}\partial \widetilde{x}_{k}}\text{.}
\end{equation*}
In terms of the $x$ variable, this is precisely the Hessian of $f,$ as
discussed above. Therefore the cycle $C\subset M$ being special is
equivalent to 
\begin{equation*}
\func{Im}\left( \widetilde{\omega }+F_{A}\right) ^{m}=\left( \tan \theta
\right) Re\left( \widetilde{\omega }+F_{A}\right) ^{m}.
\end{equation*}

For a general type-A supersymmetric cycle in $M$, we have a special
Lagrangian $C$ in $M$ together with a flat $U\left( 1\right) $ connection on
it. Since as before, $C$ is expressed as a section of $\pi :M\rightarrow B$
and is given by $y^{j}=\phi ^{jk}\frac{\partial f}{\partial x^{k}}$, a flat $%
U\left( 1\right) $ connection on $C$ can be written in the form $%
d+ide=d+i\Sigma \frac{\partial e}{\partial x^{k}}dx^{k}$ for some function $%
e=e\left( x\right) $. Recall that the transformation of $C$ alone is the
connection $d+i\Sigma y^{j}d\widetilde{y}_{j}$ over $W$. When the flat
connection on $C$ is also taken into account, then the total transformation
becomes 
\begin{eqnarray*}
D_{A} &=&d+i\Sigma y^{j}d\widetilde{y}_{j}+ide \\
&=&d+i\Sigma \phi ^{jk}\frac{\partial f}{\partial x^{k}}d\widetilde{y}%
_{j}+i\Sigma \frac{\partial e}{\partial \widetilde{x}^{j}}d\widetilde{x}_{j}.
\end{eqnarray*}
Here we have composed the function $e\left( x\right) $ with the coordinate
transformation $x=x\left( \widetilde{x}\right) .$ Notice that the added term 
$\Sigma \frac{\partial e}{\partial \widetilde{x}^{j}}d\widetilde{x}_{j}$ is
exact and therefore the curvature form of this new connection is the same as
the old one. In particular $D_{A}$ satisfies 
\begin{eqnarray*}
F_{A}^{0,2} &=&0, \\
\func{Im}e^{i\theta }\left( \omega +F\right) ^{m} &=&0,
\end{eqnarray*}
so is a supersymmetric cycle of type-B in $W$. By the same reasoning, we can
couple with $C$ a flat connection on it of \textit{any} rank and we would
still obtain a non-Abelian connection $D_{A}$ on $W$ satisfying the above
equations.

In conclusion, the transform of a type-A supersymmetric section in $M$ is a
type-B supersymmetric $2m$-cycle in $W$.

\begin{remark}
The real Fourier-Mukai transform we discussed above exchanges the symplectic
and complex aspects of the two theories.

On the A-cycle side, Donaldson and Hitchin \cite{hitchin} introduce a
symplectic form on the space of of maps $Map\left( C,M\right) $ as follows.
If $\upsilon $ is a fixed volume form on a three manifold $C$, then $%
\int_{C}ev^{\ast }\omega \wedge \upsilon $ is a symplectic form on $%
Map\left( C,M\right) $ and it equips with a Hamiltonian action by the group
of volume preserving diffeomorphisms of $C$. The zero of the corresponding
moment map is precisely the Lagrangian condition on $f\in Map\left(
C,M\right) $.

If one restricts to the infinite-dimensional complex submanifold of $%
Map\left( C,M\right) $ consisting of those $f$ which satisfy $f^{\ast
}\Omega =v$, then the symplectic quotient is the moduli space of $A$-cycles.

On the B-cycle side, we consider the pre-symplectic form $\func{Im}\left[
\int_{W}\left( \widetilde{\omega }+\mathbb{F}\right) ^{m}\right] ^{\left[ 2%
\right] }$ on the space of connections $\mathcal{A}\left( W\right) $ (see
section 4.3 or compare \cite{conan2}). This form is preserved by the group
of gauge transformations and the corresponding moment map equation is the
deformed Hermitian-Yang-Mills equations.

If one restricts to the complex submanifold of $\mathcal{A}\left( W\right) $
consisting of those connections which define a holomorphic structure on the
bundle, then the symplectic quotient is the moduli space of $B$-cycles.

Notice that the above real Fourier-Mukai transform exchanges the moment map
condition on one side to the complex condition on the other side. Such
exchanges of symplectic and complex aspects are typical in mirror symmetry.
We expect this continues to hold true in general and not just for special
Lagrangian sections in the semi-flat case.
\end{remark}

\begin{remark}
Unlike the Hermitian-Yang-Mills equation, solutions to the fully nonlinear
deformed equations may not be elliptic. On the other hand, the
special Lagrangian equation is always elliptic since its solutions are
calibrated submanifolds. Nevertheless, deformed Hermitian-Yang-Mills
connections obtained from the real Fourier-Mukai transformation as above are
always elliptic.
\end{remark}

\subsection{Transformation of a multi-section}

When $C$ is a multi-section of $\pi :M\rightarrow B,$ the situation is more
complicated. For one thing, holomorphic disks may bound $C$ -- a situation
which cannot occur in the section case (see section 4.2). Here we propose to
look at a line bundle on finite cover of $W$ as the transform B-cycle. To
begin we assume that $C$ is smooth, $\pi :C\rightarrow B$ is a branched
cover of degree $r$ (as in algebraic geometry) and $\nabla _{j}\nabla
_{k}\phi =\delta _{jk}$ for simplicity. Away from ramification locus, $C$
determines $r$ unitary connections on $W$ locally, and each satisfies the
above equation. One might be tempted to take the diagonal connection on
their direct sum, so that this $U\left( r\right) $ connection satisfies a
non-Abelian analogue of the equation. However, such a connection cannot be
defined across ramification locus because of monodromy, which can
interchange different summands of the diagonal connection.

In fact, as mentioned in \cite{mmms}, it is still open (even in physics) to
find or derive a non-Abelian analogue of the deformed Hermitian-Yang-Mills
equations via string theory (though there are some natural guesses). To
remedy this problem, we would instead construct a $U\left( 1\right) $
connection over a degree $r$ cover of $W$. This finite cover $\hat{\pi}:%
\widehat{W}\rightarrow W$ is constructed via the following Cartesian product
diagram 
\begin{equation*}
\begin{array}{ccc}
\widehat{W}=C\times _{B}W & \overset{\hat{\pi}}{\longrightarrow } & W \\ 
\downarrow &  & \downarrow \\ 
C & \overset{\pi }{\longrightarrow } & B.
\end{array}
\end{equation*}

Notice that $\widehat{W}$ is a smooth manifold because $C$ is smooth and $%
\pi :C\rightarrow B$ is a branched cover. Moreover, the construction as
before determines a smooth unitary connection over $\widehat{W}$ satisfying
the above equation (with $\omega $ replaced by $\hat{\pi}^{\ast }\omega $).

We remark that employing this construction is similar to the use of
isogenies needed to define the categorical isomorphism which proves
Kontsevich's conjecture in the case of the elliptic curve \cite{pz}. The
point is that multi-sections transforming to higher-rank can be handled via
single sections giving line bundles, by imposing functoriality after pushing
forward under finite covers.

Even though this connection over $\widehat{W}$ satisfies $F^{0,2}=0$ in a
suitable sense, $\widehat{W}$ is not a complex manifold.

\subsection{Transformation for more general cycles}

In the previous sections, we considered only sections or multi-sections on $%
M $ and obtain holomorphic bundles on $W$ which satisfied the deformed
Hermitian-Yang-Mills equation. Notice that, for a Calabi-Yau threefold $M$,
a multi-section of $\pi :M\rightarrow B$ can be characterized as special
Lagrangian cycle $C$ whose image under $\pi $ is of dimension three --
namely, the whole $B$. If the image has dimension zero, then the Lagrangian
is a torus fiber plus bundle and its dual is the point (0-brane) it
represents on the corresponding dual torus. This is the basis for the
conjecture of \cite{syz}. Here we are going to look at the other cases, that
is the dimension of the image of $C$ under $\pi $ is either (i) one or (ii)
two. For simplicity we shall only look at the flat case, namely $%
M=T^{6}=B\times F$ where both $B$ and $F$ are flat three dimensional
Lagrangian tori.

Case (i), when $\dim \pi \left( C\right) =1$. The restriction of $\pi $ to $%
C $ express $C$ as the total space of an one-parameter family of surfaces.
In fact we will see that this is a product family of an affine $T^{1}$ in $B$
with an affine $T^{2}$ in $F$.

As before we denote the coordinates of $B$ (resp. $F$) by $x^{1},x^{2},x^{3}$
(resp. $y^{1},y^{2},y^{3}$). Without loss of generality, we can assume that $%
\pi \left( C\right) $ is locally given by $x^{2}=f\left( x^{1}\right) $ and $%
x^{3}=g\left( x^{1}\right) $. Moreover the surface in $C$ over any such
point is determined by $y^{3}=h\left( x^{1},y^{1},y^{2}\right) $. In
particular, $C$ is parametrized by $x^{1},y^{1}$ and $y^{2}$ locally.

The special condition $\func{Im}\left( dz^{1}\wedge dz^{2}\wedge
dz^{3}\right) |_{C}=0$ implies that $h$ is independent of $x^{1}$. Namely
the surface family $C$ is indeed a product subfamily of $M=F\times B$. Now
the Lagrangian conditions read as follow: 
\begin{eqnarray*}
1+\frac{dg}{dx^{1}}\frac{\partial h}{\partial y^{1}} &=&0, \\
\frac{df}{dx^{1}}+\frac{dg}{dx^{1}}\frac{\partial h}{\partial y^{2}} &=&0.
\end{eqnarray*}
These imply that, 
\begin{eqnarray*}
f = x^{2} &=&ax^{1}+\alpha \\
g = x^{3} &=&bx^{1}+\beta \\
h = y^{3} &=&- \frac{1}{b}y^{1}-\frac{a}{b}y^{2}+\frac{\gamma }{b}\text{.}
\end{eqnarray*}
By analyticity of special Lagrangians, these parametrizations hold true on
the whole $C$. We can therefore express $C$ as the product $C_{B}\times
C_{F} $ with $C_{F}\subset F$ being a two torus and $C_{B}\subset B$ being a
circle. Here 
\begin{eqnarray*}
C_{B} &=&\left\{ \left( x^{1},x^{2},x^{3}\right) =\left( 1,a,b\right)
x^{1}+\left( 0,\alpha ,\beta \right) \right\} , \\
C_{F} &=&\left\{ \left( y^{1},y^{2},y^{3}\right) :\left(
y^{1},y^{2},y^{3}\right) \cdot \left( 1,a,b\right) =\gamma \right\} .
\end{eqnarray*}

Since our primary interest is when $C$ is a closed subspace of $M$, this
implies that both $a$ and $b$ are rational numbers and $C$ is a three torus
which sits in $M=T^{6}$ as a totally geodesic flat torus. To describe a
supersymmetric cycle, we also need a $U\left( 1\right) $ flat connection $%
D_{A} $ over $C$. If we parametrize $C$ by coordinate functions $x^{1},y^{2} 
$ and $y^{3}$ as above, then we have 
\begin{equation*}
D_{A}=d+i\left(\widetilde{\gamma}dx^{1}+\widetilde{\alpha}dy^{2}+ \widetilde{%
\beta}dy^{3}\right) \text{,}
\end{equation*}
for some real numbers $\widetilde{\alpha},\widetilde{\beta}$ and $\widetilde{%
\gamma}$.

Now let us define the transformation of $\left( C,L\right) $. First, the
mirror of $M={B}\times {F}$ equals $W={B}\times \widetilde{F}$ where $%
\widetilde{F}$ is the dual three torus to ${F}$. A point $\left(
x_{1},x_{2},x_{3},\widetilde{y}_{1},\widetilde{y}_{2},\widetilde{y}%
_{3}\right) =\left( x,\widetilde{y}\right) \in W$ (note $x=\widetilde{x}$
here, by flatness of the metric) lies in the mirror $\left( \widetilde{C},%
\widetilde{L}\right) $ of the SUSY cycle $\left( C,L\right) $ if and only if
the flat connection which is obtained by the restriction of $D_{A}$ to $%
x\times C_{F}$ twisted by the one form $i\sum \widetilde{y}_{j}dy^{j}$ is in
fact trivial. That is, $\widetilde{\alpha }dy^{2}+\widetilde{\beta }%
dy^{3}+\sum \widetilde{y}_{j}dy^{j}=0.$ Using the equation $\left(
y^{1},y^{2},y^{3}\right) \cdot \left( 1,a,b\right) =\gamma $, we obtain that 
\begin{eqnarray*}
\widetilde{y}_{2} &=&a\widetilde{y}_{1}-\widetilde{\alpha }, \\
\widetilde{y}_{3} &=&b\widetilde{y}_{1}-\widetilde{\beta }.
\end{eqnarray*}
Or equivalently, $\widetilde{C}=C_{B}\times C_{\widetilde{F}}$ with $C_{%
\widetilde{F}}=\left\{ \left( \widetilde{y}_{1},\widetilde{y}_{2},\widetilde{%
y}_{3}\right) =\left( 1,a,b\right) \widetilde{y}_{1}-\left( 0,\widetilde{%
\alpha },\widetilde{\beta }\right) \right\} $ and $C_{B}=\left\{ \left(
x_{1},x_{2},x_{3}\right) =\left( 1,a,b\right) x_{1}+\left( 0,\alpha ,\beta
\right) \right\} $. In particular $\widetilde{C}$ is a holomorphic curve in $%
W$. The last step would be to determine the $U\left( 1\right) $ connection $%
D_{\widetilde{A}}$ on $\widetilde{C}$. By essentially the same argument as
above and the fact that there is no transformation along the base direction,
we obtain 
\begin{equation*}
D_{\widetilde{A}}=d+i\left( \widetilde{\gamma }dx_{1}-\gamma d\widetilde{y}%
_{1}\right) .
\end{equation*}
Now the deformed Hermitian-Yang-Mills equation $\func{Im}\left( \omega +F_{%
\widetilde{A}}\right) ^{n}=0$ is equivalent to $F_{\widetilde{A}}=0,$ which
is obviously true for the above connection $D_{\widetilde{A}}$ on $%
\widetilde{C}$.

Case (ii), when $\dim \pi \left( C\right) =2.$ The restriction of $\pi $ to $%
C$ express $C$ as the total space of a two-parameter family of circles. As
before, let us write the parametrizing surface $S$ in $B$ as $x^{3}=f\left(
x^{1},x^{2}\right) $ and the one dimensional fiber over any point of it by $%
y^{2}=g\left( x^{1},x^{2},y^{1}\right) $ and $y^{3}=h\left(
x^{1},x^{2},y^{1}\right) $. Namely $C$ is parametrized by $x^{1},x^{2}$ and $%
y^{1}$ locally.

Now the Lagrangian condition would imply that each fiber is an affine circle
in $T^{3}=F$, among other things. On the other hand, the special condition $%
\func{Im}\left( dz^{1}\wedge dz^{2}\wedge dz^{3}\right) |_{C}=0$ implies
that the surface $S\subset T^{3}=B$ satisfies a Monge-Amp\`{e}re equation: 
\begin{equation*}
\det \left( \nabla ^{2}f\right) =0.
\end{equation*}

We would like to perform a transformation on $C$ which would produce a
complex surface in $W$ together with a holomorphic bundle on it with a
Hermitian-Yang-Mills connection. Notice that the deformed
Hermitian-Yang-Mills equation in complex dimension two is the same as the
Hermitian-Yang-Mills equation. To transform $C$ in this case is not as
straight forward as before because the equation governing the family of
affine circles is more complicated. However the above equation should imply
that $f$ is an affine function which would then simplify the situation a lot.

\section{Correspondence of Moduli Spaces}

Vafa \cite{vafa} has argued that the topological open string theory
describing strings ending on an A-cycle is equivalent to the topological
closed-string model on a Calabi-Yau with a bundle.\footnote{%
In order to get a bundle, we consider only sections or multi-sections.}
Equating the effective string-field theories leads to the conjecture that
the ordinary Chern-Simons theory on an A-cycle be equivalent to the
holomorphic Chern-Simons theory on the transform B-cycle. Gopakumar and Vafa
have verified equality of the partition functions for the dual resolutions
of the conifold \cite{gopavafa}. Further, all structures on the moduli
spaces of branes must be equivalent.

In the previous section, we used the Fourier-Mukai transform in the
semi-flat case to identify moduli spaces of A- and B-cycles as a set. In
this section, we extend our analysis to the Chern-Simons functional for
connections which are not necessarily flat or integrable. Then we study and
relate various geometric objects on the moduli spaces related by the
transform. In particular, we verify Vafa's conjecture in this case.

\subsection{Chern-Simons functionals}

In this section, we show the equivalence of the relevant Chern-Simons
functionals for corresponding pairs of supersymmetric cycles in our
semi-flat case. In fact, we will do this ``off-shell,'' meaning that the
equivalence holds even for connections which do not satisfy the flatness or
integrability conditions, respectively. The argument is essentially the one
given on pp. 4-5 of \cite{gvlargen}.

So, instead of flat connections, we consider a general $U\left( 1\right) $
connection $d+A$ on the special Lagrangian section $C$ in $M$, then the
above transform will still produce a connection on $W$ which might no longer
be integrable. In real dimension three, flat connections of any rank can be
characterized as those connections which are critical points of the
Chern-Simons functional, 
\begin{equation*}
CS\left( A\right) =\int_{C}Tr\left( AdA+\frac{2}{3}A^{3}\right) \text{.}
\end{equation*}
To be precise, one would need to impose boundary condition or growth
condition for $A$ because $C$ is not a closed manifold.

There is also a complexified version of Chern-Simons for any holomorphic
bundle $E$ on a Calabi-Yau threefold $W$ with holomorphic three form $%
\widetilde{\Omega }$ (\cite{DT} \cite{witten}). Namely, if $\mathcal{A}$ is
a Hermitian connection on $E$ which might not be integrable, then the
holomorphic Chern-Simons functional is given by 
\begin{equation*}
CS_{hol}\left( \mathcal{A}\right) =\int_{W}Tr\widetilde{\Omega }\wedge
\left( \mathcal{A}\bar{\partial}\mathcal{A}+\frac{2}{3}\left( \mathcal{A}%
\right) ^{3}\right) \text{.}
\end{equation*}
Notice that $CS_{hol}\left( \mathcal{A}\right) $ depends only on the $\left(
0,1\right) $ component of $\mathcal{A}$. As in the real case, $\mathcal{A}$
is a critical point for the holomorphic Chern-Simons if and only if $F_{%
\mathcal{A}}^{0,2}=0$, that is an integrable connection.

As argued in \cite{vafa}\cite{witten}, the holomorphic Chern-Simons theory
on $W$ is conjectured to be \textit{mirror} to the usual Chern-Simons theory
on $C\subset M,$ with instanton corrections given by holomorphic disks on $M$
with boundary lying on $C$ (as we will see in the next section, there are no
such instantons in our setting). In fact, we can directly check that the
Fourier-Mukai transform not only sends flat connections on $C$ to integrable
connections on $W$, but it preserves the Chern-Simons functional for an
arbitrary connection on $C$ which is not necessarily flat. Moreover, this
holds true for connections of any rank over $C$.

We consider $C=\left\{ y^{j}=\phi ^{jk}\frac{\partial f}{\partial x^{k}}%
\right\} \subset M$ as in section 3.1, and now $d+A=d+ie_{k}\left( x\right)
dx^{k}$ is an arbitrary rank $r$ unitary connection on $C$. The real
Fourier-Mukai transform of $C$ alone (resp. $C$ with the above connection)
is the connection $\mathcal{A}_{0}=d+i\phi ^{jk}\frac{\partial f}{\partial
x^{k}}d\widetilde{y_{j}}$ (resp. $\mathcal{A}=d+i\left( \phi ^{jk}\frac{%
\partial f}{\partial x^{k}}d\widetilde{y_{j}}+e_{k}\phi ^{jk}d\widetilde{%
x_{j}}\right) $) on $W$.

Recall that $\left( \mathcal{A}_{0}\right) ^{0,1}$ determines a holomorphic
structure on a bundle over $W,$ which we use as the background $\bar{\partial%
}$ operator in defining the holomorphic Chern-Simons functional. That is, $%
CS_{hol}\left( \mathcal{A}\right) =CS_{hol}\left( \mathcal{A},\mathcal{A}%
_{0}\right) $. In fact if we vary $\bar{\partial}$ continuously among
holomorphic bundles, the holomorphic Chern-Simons functional remains the
same.

Now, 
\begin{equation*}
CS_{hol}\left( \mathcal{A},\mathcal{A}_{0}\right) =\int_{W}Tr\widetilde{%
\Omega }\wedge \left( \mathcal{B}\left( \bar{\partial}-\frac{1}{2}\phi ^{jk}%
\frac{\partial f}{\partial x^{k}}d\overline{\widetilde{z_{j}}}\right) 
\mathcal{B}+\frac{2}{3}\left( \mathcal{B}\right) ^{3}\right) ,
\end{equation*}
where $\mathcal{B}=\left( \mathcal{A}-\mathcal{A}_{0}\right) ^{0,1}=\frac{i}{%
2}e_{k}\phi ^{jk}d\overline{\widetilde{z_{j}}}$. Now 
\begin{eqnarray*}
\int_{W}Tr\widetilde{\Omega }\mathcal{B}\phi ^{jk}\frac{\partial f}{\partial
x^{k}}d\overline{\widetilde{z_{j}}}\mathcal{B} &=&\left( const\right)
\int_{W}Tr\widetilde{\Omega }\mathcal{B}\frac{\partial f}{\partial 
\widetilde{x}^{j}}d\overline{\widetilde{z_{j}}}\mathcal{B} \\
&=&-\left( const\right) \int Tr\left( \varepsilon ^{2}\right) dfd\widetilde{y%
}_{1}d\widetilde{y}_{2}d\widetilde{y}_{3} \\
&=&0.
\end{eqnarray*}
Here $\varepsilon =\frac{i}{2}e_{k}\phi ^{jk}d\widetilde{x}_{k}$ is a
matrix-valued one form, and therefore $Tr\left( \varepsilon ^{2}\right) =0$.
Using the fact that $\bar{\partial}\mathcal{B=}\frac{i}{2}\frac{\partial }{%
\partial \widetilde{x}^{l}}\left( e_{k}\phi ^{jk}\right) d\overline{%
\widetilde{z_{l}}}d\overline{\widetilde{z_{j}}}$, we have 
\begin{eqnarray*}
CS_{hol}\left( \mathcal{A},\mathcal{A}_{0}\right) &=&\left( const\right)
\int_{W}Tr\widetilde{\Omega }\wedge \left( \mathcal{B}\bar{\partial}\mathcal{%
B}+\frac{2}{3}\left( \mathcal{B}\right) ^{3}\right) \\
&=&\left( const\right) \int_{W}Tr\left( \varepsilon d\varepsilon +\frac{2}{3}%
\varepsilon ^{3}\right) d\widetilde{y}_{1}d\widetilde{y}_{2}d\widetilde{y}%
_{3} \\
&=&\left( const\right) \int_{x}Tr\left( AdA+\frac{2}{3}A^{3}\right) \\
&=&\left( const\right) CS\left( A\right) .\quad \square
\end{eqnarray*}

When the dimension of $M$ is odd but bigger than three, the Fourier-Mukai
transform still preserves the (holomorphic) Chern-Simons functional even
though their critical points are no longer flat (or integrable) connection.
Instead the Euler-Lagrange equation is $\left( F_{A}\right) ^{n}=0$ (or $%
\left( F_{\mathcal{A}}^{0,2}\right) ^{n}=0$) where $\dim _{\mathbb{C}%
}M=m=2n+1$.

\subsection{Graded tangent spaces}

Transforming A-cycles to B-cycles is only the first step in understanding
mirror symmetry with branes \cite{vafa}. The next step would be to analyze
the correspondence between the moduli spaces of cycles (branes). In the next
two sections, we identify the graded tangent spaces and the holomorphic $m$
forms on the two moduli spaces of supersymmetric cycles. Generally, this
would involve holomorphic disk instanton contributions for the A-cycles
(analgously to the usual mirror symmetry A-model), but in our simplified
setting we now show these are absent.\footnote{%
The following argument is a variation of the one given on pp. 25-26 of \cite
{witten}.}

Let $D$ be a holomorphic disk whose boundary lies in the special Lagrangian
section $C$. Since we are in the local case, $C$ is homeomorphic to a ball
and we can find a closed disk $D^{\prime }\subset C$ with $\partial
D=\partial D^{\prime }$. Now using the assumption that $C$ is a section and $%
\pi _{2}\left( T^{m}\right) =0$, the closed surface $D\cup D^{\prime }$ is
contractible in $M$. Therefore $\int_{D\cup D^{\prime }}\omega =0$ by Stokes
theorem. Now $\int_{D^{\prime }}\omega =0$ because $D^{\prime }$ lies inside
a Lagrangian and $\int_{D}\omega >0$ because it is the area of $D$. This is
a contradiction. Hence there are no such holomorphic disks on $M$.

First we discuss the graded tangent spaces of the moduli of A- and B-
cycles. After that we verify that the real Fourier-Mukai transform does
preserve them in the semi-flat case.

For the A side, the tangent space to the moduli of special Lagrangians can
be identified with the space of closed and co-closed one forms (we called
such forms harmonic). This is proved by McLean \cite{mclean}. We denote it
by $H^{1}\left( C,\mathbb{R}\right) $. If $D_{A}$ is a flat $U\left(
r\right) $ connection on a bundle $E$ over $C$, then the tangent space of
the moduli of such connections at $D_{A}$ can be identified with the space
of harmonic one forms with valued in $ad\left( E\right) $. We denote it by $%
H^{1}\left( C,ad\left( E\right) \right) $. When $r=1$ the spaces $%
H^{1}\left( C,ad\left( E\right) \right) $ and $iH^{1}\left( C,\mathbb{R}%
\right) $ are the same. For $r$ bigger than one, it is expected that $%
H^{1}\left( C,ad\left( E\right) \right) $ is the tangent space at the
non-reduced point $rC$: If there is a family of special Lagrangians in $M$
converging to $C$ with multiplicity $r$, then it should determine a flat $%
U\left( r\right) $ connection on an open dense set in $C$. This connection
would extend to the whole $C$ if those special Lagrangians in the family are
branched covers of $C$ in $T^{\ast }C$. Then the tangent of this moduli
space at $rC$ should be $H^{1}\left( C,ad\left( E\right) \right) $. It is
useful to verify this statement.

The graded tangent spaces are defined to be $\oplus _{k}H^{k}\left( C,%
\mathbb{R}\right) \otimes \mathbb{C}$, or more generally $\oplus
_{k}H^{k}\left( C,ad\left( E\right) \right) \otimes \mathbb{C}$, the space
of harmonic $k$ forms and those with coefficient in $ad\left( E\right) $.

On the B side, a cycle is a $U\left( r\right) $ connection $\mathcal{D}_{A}$
on a bundle $\mathcal{E}$ over $W$ whose curvature $\mathcal{F}_{A}$
satisfies $\mathcal{F}_{A}^{0,2}=0$ and $\func{Im}e^{i\theta }\left( 
\widetilde{\omega }+\mathcal{F}_{A}\right) ^{m}=0$. If we replace the
deformed Hermitian-Yang-Mills equation by the non-deformed one, then a
tangent vector to this moduli space can be identified with an element $%
\mathcal{B}$ in $\Omega ^{0,1}\left( W,ad\left( \mathcal{E}\right) \right) $
satisfying $\bar{\partial}\mathcal{B}=0$ and $\widetilde{\omega }%
^{m-1}\wedge \partial \mathcal{B}=0$. The second equation is equivalent to $%
\bar{\partial}^{\ast }\mathcal{B}=0$. That is $\mathcal{B}$ is a $\bar{%
\partial}$-harmonic form of type $\left( 0,1\right) $ on $W$ with valued in $%
ad\left( \mathcal{E}\right) $. The space of such $\mathcal{B}$ equals the
sheaf cohomology $H^{1}\left( W,End\left( E\right) \right) $ by Dolbeault
theorem, provided $W$ is compact.

It is not difficult to see that a tangent vector $\mathcal{B}$ to the moduli
of B-cycles is a \textit{deformed }$\bar{\partial}$-\textit{harmonic form}
in the following sense:\footnote{%
If the rank of $E$ is bigger than one, then we need to symmetrize the
product in the second equation, as done in \cite{conan2}.} 
\begin{eqnarray*}
\bar{\partial}\mathcal{B} &=&0, \\
\func{Im}e^{i\theta }\left( \widetilde{\omega }+\mathcal{F}_{A}\right)
^{m-1}\wedge \partial \mathcal{B} &=&0.
\end{eqnarray*}
In general a differential form $\mathcal{B}$ of type $\left( 0,q\right) $ is
called a deformed $\bar{\partial}$-harmonic form (compare \cite{conan2}) if
it satisfies 
\begin{eqnarray*}
\bar{\partial}\mathcal{B} &=&0, \\
\func{Im}e^{i\theta }\left( \widetilde{\omega }+\mathcal{F}_{A}\right)
^{m-q}\wedge \partial \mathcal{B} &=&0.
\end{eqnarray*}
We denote this space as $\widetilde{H}^{q}\left( W,End\left( \mathcal{E}%
\right) \right) $. When the connection $\mathcal{D}_{A}$ and the phase angle 
$\theta $ are both trivial, a deformed $\bar{\partial}$-harmonic form is
just an ordinary $\bar{\partial}$-harmonic form. It is useful to know if
there is always a unique deformed $\bar{\partial}$-harmonic representative
for each coholomology class in $H^{q}\left( W,End\left( \mathcal{E}\right)
\right) $. One might want to require that $\widetilde{\omega }+\mathcal{F}%
_{A}$ is positive to ensure ellipticity of the equation.

As argued in \cite{vafa}, mirror symmetry with bundles leads to an
identification between $H^{q}\left( C,ad\left( E\right) \right) \otimes 
\mathbb{C}$ and $\widetilde{H}^{q}\left( W,End\left( \mathcal{E}\right)
\right)$ for each $q.$\footnote{%
In \cite{vafa}, the author uses $H^{k}\left( W,End\left( \mathcal{E}\right)
\right) $ instead of $\widetilde{H}^{k}\left( W,End\left( \mathcal{E}\right)
\right) $.} This can be verified in our situation as follows. For simplicity
we assume that the phase angle $\theta $ is zero and $E$ is a line bundle.

First we need to define the transformation from a degree $q$ form on $%
C\subset M$ to one on $W$. We need one transformation $\Phi $ corresponding
to deformations of special Lagrangians and another $\Psi $ corresponding to
deformations of its flat unitary bundle. 
\begin{eqnarray*}
\Omega ^{q}\left( C,ad\left( E\right) \right) \otimes \mathbb{C}
&\rightarrow &\Omega ^{0,q}\left( W,End\left( \mathcal{E}\right) \right)  \\
B_{1}+iB_{2} &\mapsto &\Phi \left( B_{1}\right) +i\Psi \left( B_{2}\right) 
\end{eqnarray*}

If $B$ is a $q$-form on the section $C\subset M$, in the coordinate system
of $x$'s we write $B=\Sigma b_{j_{1}...j_{q}}\left( x\right)
dx^{j_{1}}\cdots dx^{j_{q}}$. It suffices to define the transformations for $%
B=dx^{j},$ by naturality. The first (resp. second) transformation of $dx^{j}$
is given by $\Phi \left( B\right) =\left( \phi ^{jk}d\widetilde{y}%
_{k}\right) ^{0,1}$ (resp. $\Psi \left( B\right) =\left( dx^{j}\right)
^{0,1}=\left( \phi ^{jk}d\widetilde{x}_{k}\right) ^{0,1}$). When $q$ equals
one, these transformations are compatible with our identification of moduli
spaces of A- and B-cycles.

Now let $B=\Sigma b_{j}\left( x\right) dx^{j}$ be any one-form on $C\subset
M $. Then we have $\mathcal{B}=\Phi \left( B\right) =\left( \Sigma b_{j}\phi
^{jk}d\widetilde{y}_{k}\right) ^{0,1}=\frac{i}{2}\Sigma b_{j}\left( 
\widetilde{x}\right) \phi ^{jk}d\overline{\widetilde{z_{k}}}$. Therefore $%
\bar{\partial}\mathcal{B}=\frac{i}{2}\Sigma \frac{\partial }{\partial 
\widetilde{x}_{l}}\left( b_{j}\left( \widetilde{x}\right) \phi ^{jk}\right) d%
\overline{\widetilde{z_{l}}}d\overline{\widetilde{z_{k}}}$ and its vanishing
is clearly equivalent to $dB=0$ under the coordinate change $\frac{\partial 
}{\partial \widetilde{x}_{l}}=\phi ^{lk}\frac{\partial }{\partial x^{k}}$.
It is also easy to see that this equivalence between $dB=0$ and $\bar{%
\partial}\mathcal{B}=0$ holds true for any degree $q$ form too.

Our main task is to show that $d^{\ast }B=0$ if and only if $\func{Im}\left( 
\widetilde{\omega }+\mathcal{F}_{A}\right) ^{m-q}\wedge \partial \mathcal{B}%
=0$ for any degree $q$ form $B$ on $C\subset M$ and $\mathcal{B}=\Phi \left(
B\right) $. Notice that, by type considerations, the latter condition is the
same as $\func{Im}d\left[ \left( \widetilde{\omega }+\mathcal{F}_{A}\right)
^{m-q}\mathcal{B}\right] =0$.

Let $B$ be any degree $q$ form on the special Lagrangian $C$. Using the
symplectic form $\omega $, we obtain a $q$-vector field $v_{B},$ i.e. a
section of $\Lambda ^{q}\left( T_{M}\right) $. Then using arguments as in 
\cite{mclean} or \cite{hitchin}, we have 
\begin{equation*}
\ast B=\pm \iota _{v_{B}}\func{Im}\Omega \text{.}
\end{equation*}
Therefore $d^{\ast }B=0$ if and only if $\func{Im}d\left( \iota
_{v_{B}}\Omega \right) =0$ on $C$.\ If we write $B=\Sigma
b_{j_{1}...j_{q}}\left( x\right) dx^{j_{1}}\cdots dx^{j_{q}}$ then $%
v_{B}=\Sigma b_{j_{1}...j_{q}}\phi ^{j_{1}k_{1}}\cdots \phi ^{j_{q}k_{q}}%
\frac{\partial }{\partial y^{k_{1}}}\cdots \frac{\partial }{\partial
y^{k_{q}}}$. So on $C$, we have 
\begin{eqnarray*}
\iota _{v_{B}}\Omega &=&\iota _{v_{B}}\left( dx^{1}+idy^{1}\right) \wedge
\cdots \wedge \left( dx^{m}+idy^{m}\right) \\
&=&\pm \Sigma b_{j_{1}...j_{q}}\phi ^{j_{1}k_{1}}\cdots \phi
^{j_{q}k_{q}}i^{q}\prod_{l\neq k_{i}}\left( dx^{l}+idy^{l}\right) \\
&=&\pm \Sigma b_{j_{1}...j_{q}}\phi ^{j_{1}k_{1}}\cdots \phi
^{j_{q}k_{q}}i^{q}\prod_{l\neq k_{i}}\left( dx^{l}+i\frac{\partial }{%
\partial x^{p}}\left( \phi ^{lk}\frac{\partial f}{\partial x^{k}}\right)
dx^{p}\right) .
\end{eqnarray*}

Now we consider the corresponding $\mathcal{B}=\Phi \left( B\right) $ over $%
W $. Explicitly we have 
\begin{equation*}
\mathcal{B}=i^{q}\Sigma b_{j_{1}...j_{q}}\phi ^{j_{1}k_{1}}\cdots \phi
^{j_{q}k_{q}}d\overline{\widetilde{z}}_{_{k_{1}}}\cdots d\overline{%
\widetilde{z}}_{_{k_{q}}}.
\end{equation*}
Therefore we consider the form $\left( \widetilde{\omega }+\mathcal{F}%
_{A}\right) ^{m-q}\mathcal{B}$ of type $\left( m-q,m\right) $ on $W$: 
\begin{eqnarray*}
&&\left( \widetilde{\omega }+\mathcal{F}_{A}\right) ^{m-q}\mathcal{B} \\
&=&\left( \Sigma \left( \phi ^{jk}+i\frac{\partial ^{2}f}{\partial 
\widetilde{x}_{j}\partial \widetilde{x}_{k}}\right) d\widetilde{z}_{j}\wedge
d\overline{\widetilde{z}}_{k}\right) ^{m-q}\left( i^{q}\Sigma
b_{j_{1}...j_{q}}\phi ^{j_{1}k_{1}}\cdots \phi ^{j_{q}k_{q}}d\overline{%
\widetilde{z}}_{_{k_{1}}}\cdots d\overline{\widetilde{z}}_{_{k_{q}}}\right) .
\end{eqnarray*}
After the coordinate transformation $\frac{\partial }{\partial \widetilde{x}%
_{l}}=\phi ^{lk}\frac{\partial }{\partial x^{k}}$, it is now easy to see
that $\func{Im}d\left[ \left( \widetilde{\omega }+\mathcal{F}_{A}\right)
^{m-q}\mathcal{B}\right] =0$ if and only if $\func{Im}d\left( \iota
_{v_{B}}\Omega \right) =0$. That is, $B$ is a co-closed form on $C$ of
degree $q$. So $\Phi $ carries a harmonic $q$-form on $C$ to a deformed
harmonic $\left( 0,q\right) $-form on $W$. Now if $E$ is a higher rank
vector bundle over $C$, then the only changes we need in the proof are $B$
and $v_{B}$ now have valued in $ad\left( E\right) $ and we would replace the
exterior differentiation by covariant differentiation and also we need to
symmetrize the product in the $\bar{\partial}$-harmonic equation. The proof
of the equivalence for $\Psi $ is similar, and we omit it.

Therefore, we have proved that $B_{1}+iB_{2}$ is a harmonic form of degree $q
$ over $C$ if and only if $\Phi \left( B_{1}\right) +i\Psi \left(
B_{2}\right) $ is a $\bar{\partial}$-harmonic form of degree $\left(
0,q\right) $ over $W$. In particular $\Phi +i\Psi $ maps $H^{q}\left(
C,ad\left( E\right) \right) \otimes \mathbb{C}$ to $\widetilde{H}^{q}\left(
W,End\left( \mathcal{E}\right) \right) $.

\subsection{Holomorphic $m$-forms on moduli spaces}

The moduli spaces of A-cycles and B-cycles on a Calabi-Yau $m$-fold have
natural holomorphic $m$-forms. As explained in \cite{vafa}, these $m$-forms,
which inclde holomorphic disk instanton corrections on the A-cycle side, can
be identified with physical correlation functions derived from the
Chern-Simons partition function. Under Vafa's version of the mirror
conjecture with bundles, these partition functions and correlators should be
the same for any mirror pair $M$ and $W$, at least in dimension three. In
this section we recall the definitions of the holomorphic $m$-forms and
verify this equality our semi-flat case.

First we define a degree-$m$ closed form on $Map\left( C,M\right) $ by $%
\int_{C}ev^{\ast }\omega ^{m}$ where $\omega $ is the K\"{a}hler form on $M$
and $C\times Map\left( C,M\right) \overset{ev}{\longrightarrow }M$ is the
evaluation map. For simplicity we will pretend $C$ is a closed manifold,
otherwise suitable boundary condition is required. If $v$ is a normal vector
field along a Lagrangian immersion $f\in Map\left( C,M\right) $, then $v$
determines an one form $\eta _{v}$ on $C\,$. At $f\in Map\left( C,M\right) $
we have $\int_{C}ev^{\ast }\omega ^{m}\left( v_{1},...,v_{m}\right)
=\int_{C}\eta _{v_{1}}\wedge ...\wedge \eta _{v_{m}}$...

Next we need to incorporate flat connections on $C$ into the picture. We
denote $\mathcal{A}\left( C\right) $ the affine space of connections on $C$.
On $C\times \mathcal{A}\left( C\right) $ there is a naturally defined
universal connection $\mathbb{D}$ and curvature $\mathbb{F}$ (see for
example \cite{conan2}). With respect to the decomposition of two forms on $%
C\times \mathcal{A}\left( C\right) $ as $\Omega ^{2}\left( C\right) +\Omega
^{1}\left( C\right) \otimes \Omega ^{1}\left( \mathcal{A}\right) +\Omega
^{2}\left( \mathcal{A}\right) ,$ we write $\mathbb{F=F}^{2,0}+\mathbb{F}%
^{1,1}+\mathbb{F}^{0,2}$. Then for $\left( x,A\right) \in C\times \mathcal{A}%
\left( C\right) $ and $u\in T_{x}C$, $B\in \Omega ^{1}\left( C,End\left(
E\right) \right) ,$ we have $\mathbb{F}^{2,0}\left( x,A\right) =F_{A}$, $%
\mathbb{F}^{1,1}\left( x,A\right) \left( u,B\right) =B\left( u\right) $ and $%
\mathbb{F}^{0,2}=0$. We consider the following complex valued closed $m$
form on $Map\left( C,M\right) \times \mathcal{A}\left( C\right) $, 
\begin{equation*}
_{A}\Omega =\int_{C}Tr\left( ev^{\ast }\omega +\mathbb{F}\right) ^{m}\text{.}
\end{equation*}
Since the tangent spaces of the moduli of special Lagrangian and the moduli
of flat $U\left( 1\right) $ connections can both be identified with the
space of harmonic one forms\footnote{%
If the rank is greater than one, these harmonic forms will take values in
the corresponding local system.}. A tangent vector of the moduli space $_{A}%
\mathcal{M}\left( M\right) $ is a complex harmonic one form $\eta +i\mu $.
Then $_{A}\Omega $ is given explicitly as follows 
\begin{equation*}
_{A}\Omega \left( \eta _{1}+i\mu _{1},...,\eta _{m}+i\mu _{m}\right)
=\int_{C}Tr\left( \eta _{1}+i\mu _{1}\right) \wedge ...\wedge \left( \eta
_{m}+i\mu _{m}\right) .
\end{equation*}

On the $W$ side we have universal connection and curvature on the space of
connections $\mathcal{A}\left( W\right) $ as before and we have the
following complex valued closed $m$ form on $\mathcal{A}\left( W\right) $, 
\begin{equation*}
_{B}\Omega =\int_{W}\widetilde{\Omega }\wedge Tr\mathbb{F}^{m}\text{.}
\end{equation*}
As before, $_{B}\Omega $ descends to a closed $m$ form on $_{B}\mathcal{M}%
\left( W\right) $, the moduli space of holomorphic bundles, or equivalently
B-cycles, on $W$. It is conjectured by Vafa (\cite{vafa}) that under mirror
symmetry, these two forms $_{A}\Omega $ and $_{B}\Omega $ are equivalent
after instanton correction by holomorphic disks.

In our case, where $M$ is semi-flat and $C$ is a section, there is no
holomorphic disk. Also the real Fourier-Mukai transform gives mirror cycle.
Now we can verify Vafa's conjecture in this situation. Namely $_{A}\Omega $
and $_{B}\Omega $ are preserved under the real Fourier-Mukai transform.

For a closed one form on $C$ which represents an infinitesimal variation of
a $A$-cycle in $M$, we can write it as $d\eta +id\mu $ for some function $%
\eta $ and $\mu $ in $x$ variables. Under the above real Fourier-Mukai
transform, the corresponding infinitesimal variation of the mirror $B$-cycle
is $\delta \mathcal{A}=i\left( \frac{\partial \eta }{\partial \widetilde{x}%
_{j}}d\widetilde{x_{j}}+\frac{\partial \mu }{\partial \widetilde{x}_{j}}d%
\widetilde{y_{j}}\right) $. Therefore its $\left( 0,1\right) $ component is $%
\left( \delta \mathcal{A}\right) ^{0,1}=\frac{i}{2}\left( \frac{\partial
\eta }{\partial \widetilde{x}_{j}}+i\frac{\partial \mu }{\partial \widetilde{%
x}_{j}}\right) \left( d\widetilde{x_{j}}-id\widetilde{y_{j}}\right) $. So 
\begin{eqnarray*}
_{B}\Omega \left( \delta \mathcal{A}_{1},...,\delta \mathcal{A}_{m}\right)
&=&\int_{W}\widetilde{\Omega }\wedge Tr\mathbb{F}^{m}\left( \delta \mathcal{A%
}_{1},...,\delta \mathcal{A}_{m}\right) \\
&=&\int_{W}\widetilde{\Omega }\wedge \left[ \delta \mathcal{A}_{1}\wedge
\cdots \wedge \delta \mathcal{A}_{m}\right] _{sym} \\
&=&\left( const\right) \int_{W}\Pi _{j}\left( d\eta _{j}+id\mu _{j}\right) d%
\widetilde{y_{1}}\cdots d\widetilde{y_{m}} \\
&=&\left( const\right) ^{\prime }\int_{C}\Pi _{j}\left( d\eta _{j}+id\mu
_{j}\right) \\
&=&\left( const\right) _{A}^{\prime }\Omega \left( d\eta _{1}+id\mu
_{1},...,d\eta _{m}+id\mu _{m}\right) .
\end{eqnarray*}

Hence we are done. Note that the same argument also work for higher-rank
flat unitary bundles over $C$.

\section{B-cycles in $M$ are A-cycles in $T^*M$}

We now show that a B-cycle in $M$ can also be treated as a special
Lagrangian cycle in the cotangent bundle $X=T^{\ast }M\overset{\pi }{%
\rightarrow }M$. We include this observation for its possible relevance to
the recovery of ``classical'' mirror symmetry from the version with branes,
as outlined briefly in \cite{kontsevich}. The reader is be warned that the $%
2n$-form that we use for the special condition on $X$ may not be the most
natural one.

Recall that cotangent bundle of $M$, or any manifold, carries a natural
symplectic form $\vartheta =\Sigma dx^{k}\wedge du_{k}+\Sigma dy^{k}\wedge
dv_{k}$ where $x$'s and $y$'s are local coordinates on $M$ and $u$'s and $v$%
's are the dual coordinates in $T^{\ast }M$. Moreover the conormal bundle of
submanifold $C$ in $M$ is a Lagrangian submanifold with respect to this
symplectic form on $X$. There are a couple other natural closed two-forms on 
$X$: (i) the pullback of K\"{a}hler form from $M$, namely $\pi ^{\ast
}\omega $ and (ii) the canonical holomorphic symplectic form $\vartheta
_{hol}$ via the identification between $T^{\ast }M$ and $\left( T^{\ast
}M\otimes \mathbb{C}\right) ^{1,0}$. In term of local holomorphic
coordinates $z^{1},...,z^{n}$ on $M$ we have 
\begin{eqnarray*}
\pi ^{\ast }\omega &=&i\Sigma g_{j\bar{k}}dz^{j}\wedge d\bar{z}^{k}, \\
\vartheta _{hol} &=&\Sigma dz^{j}\wedge dw_{j}\text{.}
\end{eqnarray*}
Here $z^{j}=x^{j}+iy^{j}$ and $w_{j}=u_{j}+iv_{j}$. We define the $2n$ form $%
\Theta $ on $X$ using a combination of $\pi ^{\ast }\omega $ and $\vartheta
_{hol}$: 
\begin{equation*}
\Theta =\left( \pi ^{\ast }\omega +\func{Im}\vartheta _{hol}\right) ^{n}%
\text{.}
\end{equation*}
Notice that $\Theta \bar{\Theta}=\vartheta ^{2n}$ and the restriction of $%
\Theta $ to the zero section is a constant multiple of the volume form on $M$%
...

A Lagrangian submanifold $S$ in $X$ is called special Lagrangian if the
restriction of $\Theta $ to $S$ satisfies $\func{Im}\Theta =\tan \theta 
\func{Re}\Theta $ for some phase angle $\theta $. Equivalently $\func{Im}%
e^{i\theta }\Theta $ vanishes on $S$.

Next we consider a Hermitian line bundle $L$ over $M$. Let $D_{A}$ be a
Hermitian integrable connection on $L$, that is $F_{A}^{2,0}=0$. With
respect to a holomorphic trivialization of $L$, we can write $D_{A}=\partial
+\bar{\partial}+\partial \phi $ locally for some real valued function $\phi
\left( z,\bar{z}\right) $. This determines a Lagrangian submanifold $%
S=\left\{ w_{j}=\frac{\partial \phi }{\partial z^{j}}:j=1,...,n\right\} $ in 
$X$ with respect to $\vartheta $. Notice that the definition of $S$ depends
on the holomorphic trivialization of $L$. The restriction $\vartheta _{hol}$
to $S$ equals 
\begin{eqnarray*}
\vartheta _{hol}|_{S} &=&\Sigma dz^{j}\wedge d\left( \frac{\partial \phi }{%
\partial z^{j}}\right) \\
&=&\Sigma dz^{j}\wedge \left( \frac{\partial ^{2}\phi }{\partial
z^{j}\partial z^{k}}dz^{k}+\frac{\partial ^{2}\phi }{\partial z^{j}\partial 
\bar{z}^{k}}d\bar{z}^{k}\right) \\
&=&\Sigma \frac{\partial ^{2}\phi }{\partial z^{j}\partial \bar{z}^{k}}%
dz^{j}\wedge d\bar{z}^{k}.
\end{eqnarray*}
This form is pure imaginary because $\phi $ is a real valued function.
Therefore the restriction of $\Theta $ to $S$ equals 
\begin{equation*}
\Theta |_{S}=\left[ \Sigma \left( ig_{j\bar{k}}+\frac{\partial ^{2}\phi }{%
\partial z^{j}\partial \bar{z}^{k}}\right) dz^{j}\wedge d\bar{z}^{k}\right]
^{n}=\left( \omega +F\right) ^{n}.
\end{equation*}
Therefore $S$ is a special Lagrangian in $X$ if and only if $\left(
L,D_{A}\right) $ satisfies the deformed Hermitian-Yang-Mills equation on $M$.

\vskip 0.1in

{\scriptsize Naichung Conan Leung, School of Mathematics, University of
Minnesota, Minneapolis, MN 55455. (LEUNG@MATH.UMN.EDU) }

{\scriptsize Shing-Tung Yau, Department of Mathematics, Harvard University,
Cambridge, MA 02138. 
(YAU@MATH.NWU.EDU)}

{\scriptsize Eric Zaslow, Department of Mathematics, Northwestern
University, 2033 Sheridan Road, Evanston, IL 60208.
(ZASLOW@MATH.NORTHWESTERN.EDU)}

\end{document}